\documentclass[11pt]{article}
\usepackage{amsmath}
\usepackage[colorlinks=true]{hyperref}
 \hypersetup{urlcolor=blue, citecolor=blue}

\usepackage{amssymb}

\makeatletter
\renewcommand{\@begintheorem}[2]{\begin{trivlist}
\item[\hspace{\labelsep}{\bf \mbox{~~~}#1\ #2.}]}
\renewcommand{\@opargbegintheorem}[3]{\begin{trivlist}
\item[\hspace{\labelsep}{\bf \mbox{~~~}#1\ #2 {\rm (#3).}}]}
\renewcommand{\@endtheorem}{\end{trivlist}}
\makeatother

\newtheorem{theorem}{Theorem}[section]
\newtheorem{corollary}[theorem]{Corollary}

\newtheorem{lemma}[theorem]{Lemma}

\newtheorem{definition}[theorem]{Definition}
\newtheorem{remark}[theorem]{Remark}

\newtheorem{example}[theorem]{Example}

\begin{document}

\begin{center}{\Large \noindent \bf Lyapunov spectra behavior for linear discrete-time systems \\ under small perturbations}
\end{center}

\bigskip

\centerline{\scshape Adam Czornik}
\medskip
{\footnotesize
 \centerline{Department of Automatic Control and Robotics, Silesian University of Technology}
      \centerline{44-100 Gliwice, Poland, e-mail: adam.czornik@polsl.pl}
} 

\medskip

\centerline{\scshape Evgeni\u{\i} Makarov}
\medskip
{\footnotesize
 \centerline{Institute of Mathematics, National Academy of
Sciences of Belarus}
      \centerline{220072 Minsk, Belarus, e-mail: jcm@im.bas-net.by}
} 

\medskip

\centerline{\scshape Micha{\l } Niezabitowski}
\medskip
{\footnotesize
 \centerline{Department of Automatic Control and Robotics, Silesian University of Technology}
      \centerline{44-100 Gliwice, Poland, e-mail: michal.niezabitowski@polsl.pl}
} 

\medskip

\centerline{\scshape Svetlana Popova}
\medskip
{\footnotesize
 \centerline{Udmurt State University}
      \centerline{426034 Izhevsk, Russia, e-mail: udsu.popova.sn@gmail.com}
}

\medskip

\centerline{\scshape Vasilii Zaitsev}
\medskip
{\footnotesize
 \centerline{Udmurt State University}
      \centerline{426034 Izhevsk, Russia, e-mail: verba@udm.ru}
}

\bigskip

\begin{abstract}
We investigate the behavior of the Lyapunov spectrum of a linear discrete-time system under the action 
of small perturbations in order to obtain some verifiable conditions for stability and openness of the Lyapunov spectrum. 
To this end we introduce the concepts of broken away solutions and splitted systems.
The main results obtained are a necessary condition for stability and a sufficient condition 
for the openness of the Lyapunov spectrum, which is given in terms of the system itself. 
Finally, examples of using the obtained results are presented.
\end{abstract}

{\bf MSC2020:} 39A06, 39A22, 39A30.

\medskip

{\bf Keywords:} Discrete linear time-varying system, Lyapunov spectrum, small perturbations, multiplicative perturbations.

\section{Introduction}

The concept of characteristic exponents of linear
time-varying differential equations was introduced by A.M.~Lyapunov in 1892
in his famous work~\cite{Lyapunov}. Subsequently, the theory of Lyapunov
characteristic exponents 
developed
into a well-established asymptotic
theory of linear systems \cite{Adrianova,Barreira,BVGN,Izobov,Izobovk}. The
characteristic number $\lambda$, 
or \emph{Lyapunov exponent}, as it is called nowadays, of a
nonzero function $\varphi \colon \lbrack t_{0},+\infty )\to \mathbb{R}^{s}$ characterizes its growth 
as $t$ tends to  $+\infty $
 in the scale of
exponents $\alpha $ of exponential functions $e^{\alpha t}$, where $\alpha \in\mathbb{R}$, 
and it is defined to be
\begin{equation*}
\lambda \lbrack \varphi ]=\limsup_{t\to +\infty }t^{-1}\ln \|
\varphi (t)\| .
\end{equation*}%
A.M.~Lyapunov showed that if the original differential system%
\begin{equation}
\dot{x}=f(t,x),\quad t\ge t_{0},~x\in \mathbb{R}^{s},  
\label{X8}
\end{equation}%
has the trivial solution and Lyapunov exponents of all solutions of 
the linearized system%
\begin{equation}
\dot{x}=A(t)x,\quad t\ge t_{0},~x\in \mathbb{R}^{s},  
\label{ContTime}
\end{equation}%
where $A(t)=f'_{x}(t,0)$,
are negative, then under certain conditions on function $f(t,x)$, the trivial
solution of system~\eqref{X8} is asymptotically stable. This result was a
basis of the so-called first Lyapunov method of studying 
the stability of solutions of differential systems.
 When studying Lyapunov exponents of linear
systems, some unexpected effects were discovered, which
Lyapunov himself probably did not suspect.
 In particular, it turned out that small
perturbations of the coefficient matrix $A(\cdot )$ of system~\eqref{ContTime} 
may lead to jumps of the Lyapunov exponents of its
solutions. For example, the original linear system~\eqref{ContTime} may have
all solutions with negative Lyapunov exponents (and for this reason to be
exponentially stable), but an arbitrarily small additive perturbation 
$Q(\cdot) $ of the matrix $A(\cdot )$ 
may result in positive Lyapunov exponents of the perturbed system
\begin{equation}
\dot{x}=\bigl(A(t)+Q(t)\bigr)x,\quad t\ge t_{0},~x\in \mathbb{R}^{s}.
\label{perturbed}
\end{equation}%
 This instability phenomenon of 
Lyapunov exponents was discovered by O.~Perron~\cite{Perron} in 1930. 

All of the above also applies to a linear discrete time-varying system 
\begin{equation}
x(n+1)=A(n)x(n),\quad n\in \mathbb{N},~x\in \mathbb{R}^{s},  \label{2}
\end{equation}%
with the assumption that $A(\cdot)$ is a Lyapunov sequence
(see below, in the notation section, for the definition). 
It can be proved~\cite{Ban} that there are systems of the form~\eqref{2} for which small perturbations 
of the sequence $A(\cdot)$ lead to a significant change in the Lyapunov spectrum.

The stability property of the Lyapunov spectrum of system~\eqref{ContTime} (or system~\eqref{2}) 
means that small perturbations of the matrix $A(\cdot)$ produce a small change in the Lyapunov spectrum.
Necessary and sufficient
conditions for stability of the Lyapunov spectrum
of a continuous-time system of the form~\eqref{ContTime} were obtained by V.M.~Millionshchikov~\cite{Million2} 
and at the same time by B.F.~Bylov and N.A.~Izobov~\cite{BylovIzobov} using the so-called 
Millionschikov rotation method (see~\cite{Million}). 
They showed that 
in order for the
stability of the Lyapunov spectrum of system~\eqref{ContTime} to hold, it is necessary
and sufficient that this system can be reduced to a block-triangular form by
some Lyapunov transformation, 
such that the blocks are integrally separated~\cite{Bylov} and
for each block, the upper and lower central exponents~\cite{BylovIzobov,BVGN,Million2} coincide.
Similar conditions hold for linear discrete-time systems~\cite{BP}.
It is important to notice that, in general these conditions are
unverifiable, since for their application we must transform system~\eqref{ContTime} (or system~\eqref{2}) 
into some special form by Lyapunov transformation, but the algorithms to construct this transformation are unknown. 
Therefore, the question arises: is it possible to obtain any stability conditions for the Lyapunov spectrum
that are, in some sense, verifiable?

If the Lyapunov spectrum is stable, we are sure that sufficiently small perturbations 
of the system do not remove its Lyapunov spectrum from some small neighborhood of the original spectrum.
The following question arises naturally in this connection: 
is it possible to move the Lyapunov spectrum to any prescribed position in a small vicinity of the original spectrum 
using appropriate small perturbations?
This property can be called the openness of the Lyapunov spectrum of system~\eqref{ContTime} (or system~\eqref{2}).
Some results on the openness for continuous-time systems were obtained in~\cite{MP} and for discrete-time systems
in~\cite{BP}. 

The foundations of the asymptotic theory of discrete-time systems 
including the issues close to the problems formulated above
are presented in~\cite{Czornik,Elaydi,GAI}. 
Questions on integral separateness  and stability of Lyapunov spectrum of system~\eqref{2} 
were considered in~\cite{Ban,BP,IR}. In this context,  
we also mention the paper of L.~Barreira and C.~Valls~\cite{Barreira2}, where the problem of coincidence
of Lyapunov spectra of perturbed and unperturbed systems is investigated.
Note that these results do not allow us to achieve the goals of our article, which are to study the behavior
of the Lyapunov spectrum of system~\eqref{2} under small perturbations of the matrix $A(\cdot)$. 
We investigate a necessary condition for
stability of the Lyapunov spectrum of system~\eqref{2}, which do not require
a reduction of this system to a special form, but is expressed in terms of the system itself.
In addition, we obtain sufficient conditions for the openness of the spectrum of system~\eqref{2}. 
To solve these problems, we introduce and use the concept of splitness of system~\eqref{2} based on 
the angular behavior of solutions of this system. 

This article is organized as follows. In Section \ref{sect2}, the necessary notation
and the concept of the Lyapunov spectrum of system~\eqref{2} 
is introduced.
In Section~\ref{sect3}, the concept of the spectral set of this system under various
perturbations of its coefficient matrix is considered.
It is demonstrated that
multiplicative perturbations are more adequate to the problem under consideration.
A definition
of the stability of the Lyapunov spectrum is also introduced. In Section~\ref{sect4}
the concept of splitted systems is proposed and their properties are
discussed. In Section~\ref{sect5}, the main theorem on the property of 
 splitted systems are shown.
In Section~\ref{sect6}, we prove several results that follow from 
the main theorem
and demonstrate the importance of the introduced 
concept of splitted systems for studying the behavior of the Lyapunov spectrum  under the action of 
small perturbations. 
Some examples are given in Section~\ref{sect6}. The article is completed by Conclusions.

\section{Notation}\label{sect2}

Let $\mathbb{R}^{s}$~be an $s$-dimensional Euclidean space with a fixed
orthonormal basis $e_{1},\ldots ,e_{s}$ and the standard norm $\| \cdot
\| $. By $\mathbb{R}^{s\times t}$ we shall denote the space of all real
matrices of  size $s\times t$ with the spectral norm, i.e., with the
operator norm generated in $\mathbb{R}^{s\times t}$ by Euclidean norms in $%
\mathbb{R}^{s}$ and $\mathbb{R}^{t}$, respectively. By $[ a_{1},\ldots
,a_{t}] \in \mathbb{R}^{s\times t}$ we denote the matrix with the
sequential columns $a_{1},\ldots ,a_{t}\in \mathbb{R}^{s}$; $I\in \mathbb{R}%
^{s\times s}$ is the identity matrix. For any nonsingular matrix $H\in 
\mathbb{R}^{s\times s}$  we denote by $\varkappa (H)$ the condition number of 
$H$ with respect to spectral norm, i.\thinspace e., $\varkappa (H)=\|
H\| \,\| H^{-1}\| $. For any sequence $F(\cdot )=\bigl(F(n)\bigr)%
_{n\in \mathbb{N}}$, where $F(n)\in \mathbb{R}^{s\times t}$, $n\in \mathbb{N}
$, we define $\| F\| _{\infty }=\sup\limits_{n\in \mathbb{N}}\|
F(n)\| $. Any bounded sequence $L(\cdot )$ of invertible matrices $%
L(n)\in \mathbb{R}^{s\times s}$, $n\in\mathbb{N}$, such that the sequence $%
L^{-1}(\cdot )$ is bounded on $\mathbb{N}$, is called 
\emph{a Lyapunov sequence}.

By $\mathbb{R}_{\le }^{s}$ (resp. $\mathbb{R}_{<}^{s}$) we denote the set of all
nondecreasing (resp. increasing) sequences of $s$ real numbers. For a fixed
sequence $\mu =(\mu _{1},\dots ,\mu _{s})\in \mathbb{R}_{\le }^{s}$ and any 
$\delta >0$ let us denote by $\mathcal{O}_{\delta }(\mu )$ the set of all
sequences $\nu =\bigl(\nu _{1},\dots ,\nu _{s}\bigr)\in \mathbb{R}_{\le
}^{s}$ such that $\max_{j=1,\ldots ,s}|\nu _{j}-\mu _{j}|<\delta $. In other
words, $\mathcal{O}_{\delta }(\mu )$ is a $\delta $-neighborhood of the
sequence $\mu \in \mathbb{R}_{\le }^{s}$ with respect to the metric
generated by the vector $l_{\infty }$ norm of the space $\mathbb{R}^{s}$ on
its subset $\mathbb{R}_{\le }^{s}$.

By $[\alpha ]$ we shall denote the integer part of $\alpha \in \mathbb{R}$,
that is, $[\alpha ]$ is the largest integer not exceeding~$\alpha $.

Let us define the angle between a nonzero vector $p\in \mathbb{R}^{s}$ and
some non-trivial linear subspace $V\subset \mathbb{R}^{s}$ by the 
equality
\begin{equation*}
\sphericalangle (p;V)=\inf\limits_{0\neq q\in V}\sphericalangle (p,q),
\end{equation*}%
where 
\begin{equation*}
\sphericalangle (p,q)=\arccos \frac{\langle p,q\rangle }{\| p\|
\,\| q\| }
\end{equation*}%
is the angle between the nonzero vectors $p,q\in \mathbb{R}^{s}$, $\langle
p,q\rangle $ is the scalar product of the vectors $p$ and $q$.

In our further considerations we shall use the following lemmas. 

\begin{lemma}[\protect\cite{Popova2003}]
\label{L1}
\it
Let $V$ be a linear subspace of $\mathbb{R}^{s}$, $\dim V=s-1$;
and let $p\in $ $\mathbb{R}^{s}\backslash V$. If $X\in \mathbb{R}^{s\times
s} $ is a nonsingular matrix, then 
\begin{equation*}
\sphericalangle (Xp;XV)\ge \frac{2}{\pi }\sphericalangle (p;V)(
\varkappa (X)) ^{1-s}.
\end{equation*}
\end{lemma}

\begin{lemma}[\protect\cite{Popova2003}]
\label{L3}
\it
Let $V$ be a linear subspace of $\mathbb{R}^{s}$, $\dim V=s-1$;
and let $p\in $ $\mathbb{R}^{s}\backslash V$ be an arbitrary nonzero vector.
If a matrix $H\in \mathbb{R}^{s\times s}$ satisfies the conditions $Hp=p$
and $Hx=0$ for each $x\in V$, then 
\begin{equation*}
\| H\| =\dfrac{1}{\sin \sphericalangle (p;V)}.
\end{equation*}
\end{lemma}

\begin{lemma}[\protect\cite{Popova2003}]
\label{L2}
\it
Let $a\colon \mathbb{N}\to \mathbb{R}$ and $b\colon 
\mathbb{N}\to \mathbb{R}$ be arbitrary bounded mappings, and let 
\begin{equation*}
\psi (\mu )=\limsup_{k\to \infty }\bigl(a(k)+\mu b(k)\bigr).
\end{equation*}%
Then the following assertions are valid:

$(1)$ The function $\psi\colon\mathbb{R}\to\mathbb{R}$ is convex and
satisfies the Lipschitz condition on $\mathbb{R}$.

$(2)$ If there exists a strictly increasing sequence $( k_{j})
_{j\in \mathbb{N}}$ of positive integers such that%
\begin{equation*}
\lim_{j\to \infty }a( k_{j}) =\psi ( 0) ,\
\rho \doteq \lim_{j\to \infty }b(k_{j})>0,
\end{equation*}%
then the function $\psi (\cdot )$ is (strictly) monotone increasing on the
interval $[ 0,\infty ) $ and the estimate%
\begin{equation*}
\psi ( \mu ) -\psi ( 0) \ge \rho \mu
\end{equation*}%
is valid for all $\mu \ge 0$. Moreover for each $t\ge 0,$ there exists a $%
\mu _{t},$ $0\le \mu _{t}\le \rho ^{-1}t,$ such that 
\begin{equation*}
\psi ( \mu _{t}) =\psi ( 0) +t.
\end{equation*}
\end{lemma}

Consider a discrete linear time-varying system~\eqref{2} with a Lyapunov
sequence $A\colon \mathbb{N}\to \mathbb{R}^{s\times s}$. Put 
\begin{equation*}
a\doteq \max \bigl\{\| A\| _{\infty },\| A^{-1}\| _{\infty }%
\bigr\}<\infty .
\end{equation*}%
Note that 
\begin{equation*}
\| A\| _{\infty }+\| A^{-1}\| _{\infty }\ge \| A(1)\|
+\| A^{-1}(1)\| 
\ge \| A(1)\| +\| A(1)\| ^{-1}\ge 2,
\end{equation*}%
hence $a\ge 1$.

We denote the transition matrix of system~\eqref{2} by $X_{A}(n,m)$, $n,m\in 
\mathbb{N}$. Then \cite[p.\,13--14]{GAI} for each solution $x(\cdot )$ of
this system, we have the equality%
\begin{equation*}
x(n)=X_{A}(n,m)x(m)\quad \text{for all\ }n\in \mathbb{N},\ m\in \mathbb{N},
\end{equation*}%
and 
\begin{equation*}
X_{A}(n,m)=\left\{ 
\begin{array}{ll}
\prod\limits_{l=m}^{n-1}A(l) & \text{ for\ }n>m, \\ 
I & \text{ for\ }n=m, \\ 
X_{A}^{-1}(m,n) & \text{ for\ }n<m.%
\end{array}%
\right.
\end{equation*}%
Then for any $n\in \mathbb{N}$, $m\in \mathbb{N}$ the following inequality%
\begin{equation}
\| X_{A}(n,m)\| \leqslant a^{|n-m|}  \label{otsMK}
\end{equation}%
is true.

Note that here and 
throughout the paper we put 
\begin{equation*}
{\prod\limits_{l=m}^{n-1}A(l)=A(n-1)A(n-2)\ldots A(m)},
\end{equation*}%
i.e., the matrices are multiplied in descending order of the index.

Let $\Phi (\cdot )$ be a fundamental system of solutions (FSS) of system~%
\eqref{2}, i.e., a set of $s$ linearly independent solutions $x_{1}(\cdot
),\dots ,x_{s}(\cdot )$ of system~\eqref{2}. We identify FSS $\Phi(\cdot)$ with the matrix 
$\bigl[ x_{1}(\cdot ),\dots ,x_{s}(\cdot )\bigr] $, which is called a 
fundamental matrix (FM)  of system~\eqref{2}.

For any nontrivial solution $x(\cdot )$ of system~\eqref{2} the Lyapunov
exponent $\lambda \lbrack x]$ is defined as 
\begin{equation*}
\lambda \lbrack x]=\limsup_{n\to \infty }\frac{1}{n}\ln \|
x(n)\| .
\end{equation*}%
It is well known~\cite{Barreira} that if $A(\cdot )$ is a Lyapunov sequence,
then the set of Lyapunov exponents of all nontrivial solutions of system~%
\eqref{2} are included in the 
interval $[-\ln a,\ln a]$, 
and 
contains at most $s$ 
elements, say 
\begin{equation*}
\Lambda _{1}(A)<\Lambda _{2}(A)<\ldots <\Lambda _{r}(A).
\end{equation*}%
The Lyapunov exponent of the trivial solution of system~\eqref{2} is set
equal to $-\infty $.

For each $j\in \{1,\dots ,r\}$ let us consider the set $\mathcal{E}_{j}$ of
all solutions of system \eqref{2}, whose Lyapunov exponents do not exceed $%
\Lambda _{j}$. Moreover, by $\mathcal{E}_{0}$ we denote the set that
consists of the trivial solution of system \eqref{2}.

Then \cite[p.\,54]{GAI} each of the sets $\mathcal{E}_{j}$ is a linear
subspace, and the dimension of the subspace $\mathcal{E}_{j}$ is equal to $%
N_{j}$, where $N_{j}$~is the maximal number of linearly independent solutions of
system \eqref{2}, which have Lyapunov exponents $\Lambda _{j}$. We put 
$N_{0}=\dim \mathcal{E}_{0}\doteq 0$. Since 
$\mathcal{E}_{j}\subset \mathcal{E}_{l}$ for $j<l$, then $N_{0}<N_{1}<\ldots <N_{r}$, and $N_{r}=s$.

Let $\Phi (\cdot )=\bigl\{x_{1}(\cdot ),\ldots ,x_{s}(\cdot )\bigr\}$ be an
arbitrary FSS of system~\eqref{2}. For each $j\in \{1,\dots ,r\}$ consider
the value $s_{j}$ which is the number of solutions from the set $\Phi (\cdot
)$ with exponent equal to $\Lambda _{j}$. It is known \cite[p.\,54]{GAI},
that the following inequalities 
hold:
\begin{equation*}
s_{1}+\ldots +s_{j}\le N_{j},\quad j=1,\dots ,r.
\end{equation*}%

\begin{definition}[{\protect\cite[p.\,53]{GAI}}]
\label{normal}FSS $\Phi (\cdot )$ is called \emph{normal}, if the
following equalities hold:
\begin{equation*}
s_{1}+\ldots +s_{j}=N_{j},\quad j=1,\dots ,r.
\end{equation*}%
\end{definition}

It is known~\cite[p.\,55]{GAI}, that for each system~\eqref{2} a normal FSS
exists.

\begin{definition}[{\protect\cite[p.\,55]{GAI}}]
\label{neszhim}We say that the FSS $\Phi (\cdot )=\left\{ x_{1}( \cdot
) ,\ldots ,x_{s}( \cdot ) \right\} $ is \emph{%
incompressible}, if for any nontrivial combination $y(\cdot
)=\sum\limits_{j=1}^{s}c_{j}x_{j}(\cdot )$ the
equality
\begin{equation*}
\lambda \lbrack y]=\max \bigl\{\lambda \lbrack x_{j}]\colon c_{j}\neq 0%
\bigr\}
\end{equation*}%
holds.
\end{definition}

It is known (see~\cite[p.\,55]{GAI}), that a FSS $\Phi (\cdot )$ is normal
if and only if it is incompressible.

Definition~\ref{normal} implies an important consequence: for each normal
FSS of system~\eqref{2}, the number $s_{j}$ of its solutions with the
Lyapunov exponent $\Lambda _{j}$ is the same and coincides with the value $%
N_{j}-N_{j-1}$, $j=1,\dots ,r$. Thus, we can associate with each linear
discrete  time-varying system~\eqref{2} a collection of $s$ numbers $%
\lambda _{1},\lambda _{2},\ldots ,\lambda _{s}$, which are the Lyapunov
exponents of the solutions included in any normal FSS of our system. This
collection is called \emph{the Lyapunov spectrum} of system~\eqref{2}~\cite%
[p.\,57]{GAI}. Further we denote it by 
\begin{equation*}
\lambda (A)=\bigl(\lambda _{1}(A),\ldots ,\lambda _{s}(A)\bigr),
\end{equation*}%
assuming that the inequalities $\lambda _{1}(A)\le \ldots \le \lambda
_{s}(A)$ are satisfied, and therefore $\lambda (A)\in \mathbb{R}_{\leqslant
}^{s}$.

\section{Preliminaries} \label{sect3}
Let us consider an additively perturbed system
\begin{equation}
x(n+1)=\bigl(A(n)+Q(n)\bigr)x(n),\quad n\in \mathbb{N},\ x\in \mathbb{R}^{s},
\label{4}
\end{equation}%
with $A( \cdot ) $ being the Lyapunov sequence and $Q(\cdot )$
being the additive perturbation.

\begin{definition}[\protect\cite{BP}]
\label{DopAdd}A sequence $Q(\cdot )$ is said to be \emph{an admissible
additive perturb\-ation} for system~\eqref{2} if $A(\cdot )+Q(\cdot )$ is a
Lyapunov sequence.
\end{definition}

Since $A(\cdot)$ is a Lyapunov sequence, it is easy to see that the following lemma holds.

\begin{lemma}[\protect\cite{BP}]
\label{lem1} 
\it
Sequence $Q(\cdot )$ is an admissible additive perturbation for
system~\eqref{2} if and only if there exists a Lyapunov sequence $R\colon
\mathbb{N}\to \mathbb{R}^{s\times s}$ such that
\begin{equation}
Q(n)=A(n)R(n)-A(n),\quad n\in \mathbb{N}.  \label{M2}
\end{equation}
\end{lemma}

Under the assumption that $Q(\cdot)$ is an admissible additive perturbation,
equality~\eqref{M2} enables us to rewrite the perturbed system~\eqref{4} in
the following form%
\begin{equation}
x(n+1)=A(n)R(n)x(n),\quad n\in \mathbb{N},\ x\in \mathbb{R}^{s}.  \label{M4}
\end{equation}%
Here the sequence $R(\cdot )$ is called \emph{a multiplicative perturbation%
} whereas system~\eqref{M4} is called  \emph{a multiplicatively perturbed
system}. 
Since,
according to our assumptions $A(\cdot )$ is a Lyapunov sequence then we
arrive
at the following definition.

\begin{definition}[\protect\cite{BP}]
\label{DopMult}A multiplicative perturbation $R(\cdot )$ is said to be
\emph{an admissible multiplicative perturbation} if $R(\cdot )$ is a
Lyapunov sequence.
\end{definition}

\begin{remark}
\label{zam1} Let us notice that for all systems~\eqref{2} the set of
admissible multiplicative perturbations is the same, whereas the set of
admissible additive perturbations depends on the coefficient matrix $%
A(\cdot) $ of  system~\eqref{2}.
\end{remark}

For a fixed sequence $A(\cdot )$, 
let $\mathcal{Q}$ denote
 the set of all
systems~\eqref{4} corresponding to admissible 
additive perturbations $Q(\cdot )$  for system~\eqref{2} and similarly,
let $\mathcal{R}$ denote the set of
all systems~\eqref{M4} corresponding to  admissible multiplicative
perturbations $R(\cdot )$. 
We write $Q(\cdot)\in\mathcal Q$ identifying system~\eqref{4} and the additive perturbation $Q(\cdot)$ that defines this system. 
Similarly, we use the notation $R(\cdot)\in\mathcal R$ for system~\eqref{M4} and the corresponding multiplicative perturbation $R(\cdot)$.
From Lemma~\ref{lem1} and definitions of
admissible perturbations we have%
\begin{equation}
\mathcal{Q}=\mathcal{R}.  \label{Wn1}
\end{equation}
We also use some subsets of the sets $\mathcal{Q}$ and $\mathcal{R}$. 
For any $\delta >0$ let us denote by $\mathcal{Q}_{\delta }$ the subset of 
$\mathcal{Q}$ corresponding to 
sequences $Q(\cdot )$ satisfying 
\begin{equation*}
\| Q\| _{\infty }<\delta, 
\end{equation*}%
and by $\mathcal{R}_{\delta }$ the subset of $\mathcal{R}$ corresponding to
sequences $R(\cdot )$ satisfying%
\begin{equation*}
\| R-I\| _{\infty }<\delta .
\end{equation*}

\begin{lemma}[\protect\cite{BP}]
\label{lem2}
\it
For any $\delta >0$ we have $\mathcal{Q}_{\delta }\subset
\mathcal{R}_{a\delta }$ and $\mathcal{R}_{\delta }\subset \mathcal{Q}%
_{a\delta }$, where $a=\max \left\{ \left\| A\right\| _{\infty
},\left\| A^{-1}\right\| _{\infty }\right\} $.
\end{lemma}

\emph{Proof.}
Take any $Q(\cdot)\in\mathcal{Q}_{\delta}$. Put
\begin{equation*}
R(n)=I+A^{-1}(n)Q(n),\quad n\in \mathbb{N}.
\end{equation*}
From Lemma~\ref{lem1} it follows that $R(\cdot)$ is an admissible
multiplicative perturbation for system~\eqref{2}.
Since 
\begin{equation*}
R(n)=I+A^{-1}(n)Q(n),\quad n\in \mathbb{N},
\end{equation*}
we have $R(\cdot)\in\mathcal{R}_{a\delta }$. Hence $\mathcal{Q}_{\delta }\subset
\mathcal{R}_{a\delta }$.

Now take any $R(\cdot)\in\mathcal{R}_{\delta }$. Put
\begin{equation*}
Q(n)=A(n)R(n)-A(n),\quad n\in \mathbb{N}.
\end{equation*}
By Lemma~\ref{lem1} we see that $Q(\cdot )$ is an admissible additive
perturbation for system~\eqref{2} and%
\begin{equation*}
\| Q\| _{\infty }\le \| A\| _{\infty }\| R-I\| _{\infty
}<a\delta .
\end{equation*}%
Hence $Q(\cdot)\in\mathcal{Q}_{a\delta }$ and, therefore, $\mathcal{R}_{\delta }\subset \mathcal{Q}%
_{a\delta }$.
\hfill $\square$

If $Q(\cdot)\in\mathcal Q$,
then we can define the Lyapunov spectrum $\lambda (A+Q)\in \mathbb{R}_{\le}^{s}$ of system~\eqref{4}.
In a similar way we can define the Lyapunov spectrum $\lambda (AR)\in \mathbb{R}_{\le}^{s}$ of system~\eqref{M4}
for $R(\cdot)\in\mathcal R$.

The spectral set of system~\eqref{2} corresponding to the class $\mathcal{Q}$ is defined by the equality
$\lambda(\mathcal{Q})\doteq\bigl\{\lambda(A+Q)\colon Q(\cdot)\in\mathcal Q\bigr\}$.
Similarly, we define the sets
$\lambda(\mathcal{Q}_{\delta})$, $\lambda(\mathcal{R})$ and $\lambda(\mathcal{R}_{\delta})$.

From~\eqref{Wn1} and 
Lemma~\ref{lem2} we get the following statement.

\begin{corollary}
\it
The equality $\lambda (\mathcal{Q})=\lambda (\mathcal{R})$ holds. 
Moreover,
for any $\delta >0$ 
the inclusions
\begin{equation*}
\lambda \bigl(\mathcal{Q}_{\delta }\bigr)\subset \lambda \bigl(\mathcal{R}_{a\delta }\bigr),\quad
\lambda \bigl(\mathcal{R}_{\delta }\bigr)\subset \lambda \bigl(\mathcal{Q}_{a\delta }\bigr)
\end{equation*}
hold.
\end{corollary}

\begin{definition}
\label{UstAdd}(\cite{BP,Czornik}) \textrm{The Lyapunov spectrum of system~%
\eqref{2} is called \emph{stable} if for any $\varepsilon >0$ there exists
$\delta >0$ such that the inclusion }$\mathrm{\lambda \bigl(\mathcal{Q}%
_{\delta }\bigr)\subset \mathcal{O}_{\varepsilon }\bigl(\lambda (A)\bigr)}$%
\textrm{\ is satisfied.}
\end{definition}

By Lemma~\ref{lem2} we obtain the following result.

\begin{theorem}[\protect\cite{BP}]
\it
The Lyapunov spectrum of  system~\eqref{2} is stable if and only if for
any $\varepsilon >0$ there exists $\delta >0$ such that the inclusion $%
\lambda \bigl(\mathcal{R}_{\delta }\bigr)\subset \mathcal{O}_{\varepsilon }%
\bigl(\lambda (A)\bigr)$ is satisfied.
\end{theorem}

\section{Splitted systems} \label{sect4}

Let $\bigl\{x_{1}(\cdot ),\dots ,x_{s}(\cdot )\bigr\}$ be some FSS of system~%
\eqref{2}. For any $i\in \left\{ 1,\dots ,s\right\} $ and $n\in \mathbb{N}$,
by $V_{i}(n)$ we denote the linear span of the vectors $x_{j}(n)$, $j\in
\{1,\dots ,s\}\setminus \{i\}$ and by $\varphi _{i}(n)\doteq \sphericalangle %
\bigl(x_{i}(n);V_{i}(n)\bigr)$ we denote the angle between the vector $%
x_{i}(n)$ and the subspace $V_{i}(n)$. We take an arbitrary $\sigma \in 
\mathbb{N}$. For any $\gamma \in ( 0,\frac{\pi }{2}] ,$ $k\in 
\mathbb{N}$, and $i\in \left\{ 1,\dots ,s\right\} $, we set 
\begin{equation*}
\Gamma _{i}^{\gamma }( \sigma ) \doteq \left\{ j\in \mathbb{N}%
\colon \varphi _{i}(j\sigma )\ge \gamma \right\} ,
\qquad
\Gamma _{i}^{\gamma }( k;\sigma ) \doteq \Gamma _{i}^{\gamma
}( \sigma ) \cap \left\{ 1,\dots ,k\right\} .
\end{equation*}%
Let $N_{i}^{\gamma }( k;\sigma ) $ be the number of elements of
the set $\Gamma _{i}^{\gamma }( k;\sigma ) $. Let us also
introduce the following notation%
\begin{equation*}
g_{i}^{\gamma }( k;\sigma ) \doteq \frac{N_{i}^{\gamma }(
k;\sigma ) }{k},
\qquad
f_{i}( k;\sigma ) \doteq \frac{\ln \left\| x_{i}( k\sigma
) \right\| }{k\sigma }.
\end{equation*}%
If the numbers $\gamma $ and $\sigma $ are given in advance, then the
corresponding symbols in the 
notation introduced above will be omitted.

A sequence $\bigl(n_{k}\bigr)_{k\in \mathbb{N}}$ of natural numbers strictly
increasing to $+\infty $ is referred to as 
\emph{a realizing sequence}
of a solution $x(\cdot )$ of system~\eqref{2} if 
\begin{equation*}
\lambda [x]=\lim_{k\to \infty }\frac{\ln \|x(n_{k})\|}{n_{k}}.
\end{equation*}

\begin{definition}
\label{ras4l} \textrm{We say that the solution $x_{j}(\cdot )$ occurring in
the FSS $\Phi (\cdot )=\bigl\{x_{1}(\cdot ),\dots ,x_{s}(\cdot )\bigr\}$ is 
\emph{$\sigma $-broken away} (from the remaining solutions of $\Phi (
\cdot ) $) if for a given $\sigma \in \mathbb{N}$, there exists a $%
\gamma \in \bigl(0,\frac{\pi }{2}\bigr]$ and a realizing sequence 
$\bigl(k_{n}\sigma \bigr)_{n\in \mathbb{N}}$ of the solution $x_j(\cdot)$, where 
$k_{n}\in \mathbb{N}$, such that 
\begin{equation*}
\lim_{n\to \infty }g_{j}^{\gamma }( k_{n};\sigma ) >0.
\end{equation*}%
A FSS $\Phi (\cdot )$ is said to be \emph{$\sigma $-splitted} if each of
the solutions of this FSS is $\sigma $-broken away. }
\end{definition}

Let us consider some basic properties of the 
notions introduced above.

\begin{theorem}
\label{Tw1}
\it 
If the solution $x_{i}(\cdot )$ occurring in the FSS $\Phi
(\cdot )=\bigl\{x_{1}(\cdot ),\dots ,x_{s}(\cdot )\bigr\}$ is $\sigma _{0}$%
-broken away for some $\sigma _{0}\in \mathbb{N}$, then it is $\sigma $%
-broken away for any $\sigma \in \mathbb{N}$.
\end{theorem}

\emph{Proof.}
Let $\gamma \in \left( 0,\frac{\pi }{2}\right] $, and let a strictly
increasing sequence $( k_{n}) _{n\in \mathbb{N}}$ be chosen so as
to satisfy the conditions%
\begin{equation*}
\lim_{n\to \infty }f_{i}( k_{n};\sigma _{0}) =\lambda 
[ x_{i}]
\end{equation*}%
and 
\begin{equation*}
\lim_{n\to \infty }g_{i}^{\gamma }( k_{n};\sigma _{0}) >0.
\end{equation*}%
To each $n\in \mathbb{N}$ we assign the nonnegative integer $l_{n}$ such that%
\begin{equation}
k_{n}\sigma _{0}\in [ l_{n}\sigma ,( l_{n}+1) \sigma )
,  \label{5}
\end{equation}%
i.e.,%
\begin{equation*}
k_{n}\frac{\sigma _{0}}{\sigma }-1<l_{n}\le k_{n}\frac{\sigma _{0}}{\sigma }%
.
\end{equation*}%
Since $( k_{n}) _{n\in \mathbb{N}}$ is a strictly increasing
sequence, it follows that $( l_{n}) _{n\in \mathbb{N}}$ is
nondecreasing. Moreover,%
\begin{equation*}
\lim_{n\to \infty }l_{n}=\infty .
\end{equation*}%
Since 
\begin{equation*}
0\le k_{n}\sigma _{0}-l_{n}\sigma <\sigma ,
\end{equation*}%
we have 
\begin{equation}
\lim_{n\to \infty }\frac{k_{n}\sigma _{0}}{l_{n}\sigma }=1
\label{X1}
\end{equation}%
and by \eqref{otsMK} we get%
\begin{equation}
\left\| X_{A}( l_{n}\sigma ,k_{n}\sigma _{0}) \right\| \le
a^{\sigma },\quad \left\| X_{A}( k_{n}\sigma _{0},l_{n}\sigma
) \right\| \le a^{\sigma }.  \label{X}
\end{equation}%
This implies that%
\begin{align*}
f_{i}( l_{n};\sigma ) = & \frac{\ln \left\| X_{A}(
l_{n}\sigma ,k_{n}\sigma _{0}) x_{i}( k_{n}\sigma _{0})
\right\| }{l_{n}\sigma }\\
\le & \frac{\ln ( \left\| X_{A}( l_{n}\sigma ,k_{n}\sigma _{0})
\right\| \left\| x_{i}( k_{n}\sigma _{0}) \right\|
) }{l_{n}\sigma }\\
\le  & \frac{1}{l_{n}\sigma }( \sigma \ln a+\ln \left\| x_{i}(
k_{n}\sigma _{0}) \right\| )
\end{align*}%
and, on the other hand,%
\begin{align*}
f_{i}( l_{n};\sigma ) \ge & \frac{1}{l_{n}\sigma }\ln \left(
\left\| X_{A}( k_{n}\sigma _{0},l_{n}\sigma ) \right\|
^{-1}\left\| x_{i}( k_{n}\sigma _{0}) \right\| \right) \\
\ge & \frac{1}{l_{n}\sigma }( -\sigma \ln a+\ln \left\| x_{i}(
k_{n}\sigma _{0}) \right\| ) .
\end{align*}%
Now, by taking into account the above equalities and \eqref{X1} we get%
\begin{align*}
\lim_{n\to \infty }\frac{1}{l_{n}\sigma }( \pm \sigma \ln a+\ln
\left\| x_{i}( k_{n}\sigma _{0}) \right\| )
= &\lim_{n\to \infty }\frac{1}{l_{n}\sigma }\ln \left\| x_{i}(
k_{n}\sigma _{0}) \right\| \\
= & \lim_{n\to \infty }\frac{1}{k_{n}\sigma _{0}}\ln \left\|
x_{i}( k_{n}\sigma _{0}) \right\| =\lambda [ x_{i}]
\end{align*}%
and we obtain%
\begin{equation*}
\lim_{n\to \infty }f_{i}( l_{n};\sigma ) =\lambda [
x_{i}] .
\end{equation*}

For any $n,m\in \mathbb{N}$ the linear subspace $V_{i}( n) $ of
vectors $x_{j}( n) $, $j\neq i$ can be represented in the form $%
V_{i}( n) =X_{A}( n,m) V_{i}( m) $. Let us
denote $\sigma _{1}=\max \left\{ \sigma _{0},\sigma \right\} $ and consider $%
n,m\in \mathbb{N}$ such that $\left| n-m\right| \le \sigma _{1}$.
Then by Lemma \ref{L1} we have%
\begin{equation*}
\begin{split}
\varphi _{i}(n)\doteq &\sphericalangle (x_{i}(n);V_{i}(n)) 
=  \sphericalangle (X_{A}( n,m) x_{i}(m);X_{A}( n,m) V_{i}(m)) \\
\ge & 2\sphericalangle (x_{i}(m);V_{i}(m))( \varkappa ( X_{A}(
n,m) ) ) ^{1-s}\frac{1}{\pi } \\
= & 2\varphi _{i}( m) ( \varkappa ( X_{A}( n,m) ) ) ^{1-s}\frac{1%
}{\pi }\ge c\varphi _{i}(m),
\end{split}
\end{equation*}%
where%
\begin{equation*}
c=\frac{2a^{2\sigma _{1}( 1-s) }}{\pi }
\end{equation*}%
by \eqref{X}.

Two cases are possible.\\
Case 1 ($\sigma <\sigma _{0}$). Let $p=\left[ \frac{\sigma _{0}}{\sigma }\right] $. 
We take a $j\in \mathbb{N}$ such that $k_{j}>2$. It follows from \eqref{5}
that 
\begin{equation*}
l_{j}\le k_{j}\frac{\sigma _{0}}{\sigma }<k_{j}\left( \left[ \frac{\sigma
_{0}}{\sigma }\right] +1\right) =k_{j}( p+1) ,
\end{equation*}%
i.e., 
\begin{equation*}
\frac{k_{j}}{l_{j}}>\frac{1}{p+1}.
\end{equation*}%
Let $m$ range over the set $\Gamma _{i}^{\gamma }(
k_{j}-2;\sigma_{0})$. The interval $\left[ \left( m-\frac{1}{2}\right)
\sigma _{0},\left( m+\frac{1}{2}\right) \sigma _{0}\right) $ for each $m$
contains at least $p$ multiples of $\sigma $. At all of these points 
the angle $%
\varphi _{i}$ is not less than $c\gamma$. In addition, all these points lie
on the real line on the left of%
\begin{equation*}
\left( k_{j}-2+\frac{1}{2}\right) \sigma _{0}<k_{j}\sigma _{0}-\sigma
_{0}<l_{j}\sigma +\sigma -\sigma _{0}<l_{j}\sigma ,
\end{equation*}%
i.e., to the left of $l_{j}\sigma $. Therefore all of them belong to the set
$\Gamma _{i}^{c\gamma }( l_{j};\sigma )$.
Then for the total number of elements of the set $\Gamma _{i}^{c\gamma
}( l_{j};\sigma ) $ we have the estimate%
\begin{equation*}
N_{i}^{c\gamma }( l_{j};\sigma ) \ge pN_{i}^{\gamma }(
k_{j}-2;\sigma _{0}) \ge p( N_{i}^{\gamma }( k_{j};\sigma
_{0}) -2)
\end{equation*}%
and%
\begin{equation*}
\begin{split}
g_{i}^{c\gamma }( l_{j};\sigma ) = & \frac{N_{i}^{c\gamma }(
l_{j};\sigma ) }{l_{j}}\ge \frac{p( N_{i}^{\gamma }(
k_{j};\sigma _{0}) -2) }{l_{j}} \\
= & \frac{pN_{i}^{\gamma }( k_{j};\sigma _{0}) }{k_{j}}\frac{k_{j}}{%
l_{j}}-\frac{2p}{l_{j}}>\frac{pg_{i}^{\gamma }( k_{j};\sigma
_{0}) }{p+1}-\frac{2p}{l_{j}}.
\end{split}
\end{equation*}%
Consequently%
\begin{equation*}
\limsup_{j\to \infty }g_{i}^{c\gamma }( l_{j};\sigma )
\ge \lim_{j\to \infty }\left( \frac{pg_{i}^{\gamma }(
k_{j};\sigma _{0}) }{p+1}-\frac{2p}{l_{j}}\right) 
=\frac{p}{p+1}\lim_{j\to \infty }g_{i}^{\gamma }( k_{j};\sigma
_{0}) >0.
\end{equation*}
Case 2 ($\sigma \ge \sigma _{0}$).  Denote $q=\left[ \frac{\sigma }{\sigma _{0}}%
\right] $. Let $j\in \mathbb{N}$ satisfy the condition $l_{j}>1$. It follows
from \eqref{5} that%
\begin{equation*}
\frac{k_{j}}{l_{j}}\ge \frac{\sigma }{\sigma _{0}}\ge \left[ \frac{\sigma 
}{\sigma _{0}}\right] =q.
\end{equation*}%
We choose some positive integer $l<l_{j}$. Then%
\begin{equation*}
\left( l+\frac{1}{2}\right) \sigma <( l+1) \sigma \le
l_{j}\sigma \le k_{j}\sigma _{0}.
\end{equation*}%
The interval $[ ( l-\frac{1}{2}) \sigma ,( l+\frac{1}{2}%
) \sigma ) $ contains at most $q+1$ multiples of $\sigma _{0}$.
If at least one of them belongs to $\Gamma _{i}^{\gamma }( k_{j};\sigma
_{0}) $, then $l\in \Gamma _{i}^{c\gamma }( l_{j};\sigma ) $. 
Consequently%
\begin{equation*}
N_{i}^{\gamma }( k_{j};\sigma _{0}) \le ( q+1)
N_{i}^{c\gamma }( l_{j};\sigma )
\end{equation*}%
and%
\begin{equation*}
g_{i}^{c\gamma }( l_{j};\sigma ) =\frac{N_{i}^{c\gamma }(
l_{j};\sigma ) }{l_{j}}\ge \frac{N_{i}^{\gamma }( k_{j};\sigma
_{0}) }{( q+1) l_{j}} 
=\frac{N_{i}^{\gamma }( k_{j};\sigma _{0}) }{k_{j}}\frac{k_{j}}{%
( q+1) l_{j}}\ge \frac{g_{i}^{\gamma }( k_{j};\sigma
_{0}) q}{q+1}.
\end{equation*}%
Therefore%
\begin{equation*}
\limsup_{j\to\infty }g_{i}^{c\gamma }( l_{j};\sigma ) \ge \frac{q%
}{q+1}\lim_{j\to\infty }g_{i}^{\gamma }( k_{j};\sigma _{0}) >0.
\end{equation*}

Now from the sequence $( l_{j}) _{j\in \mathbb{N}}$ we extract a
strictly increasing subsequence $( l_{j_{m}}) _{m\in \mathbb{N}}$
on which the limit%
\begin{equation*}
\limsup_{j\to \infty }g_{i}^{c\gamma }( l_{j};\sigma )
\end{equation*}%
is realized. The subsequence satisfies the relations%
\begin{equation*}
\lim_{m\to \infty }f_{i}( l_{j_{m}};\sigma ) =\lambda 
[ x_{i}] ,
\qquad 
\lim_{m\to \infty }g_{i}^{c\gamma }( l_{j_{m}};\sigma )
>0,
\end{equation*}%
it means that the solution $x_{i}(\cdot)$ is $\sigma$-broken away. The proof
is completed.
\hfill $\square$

\begin{corollary}
\it
\label{Cor1}If a FSS is $\sigma_{0}$-splitted for certain $\sigma_{0}\in 
\mathbb{N}$, then it is $\sigma$-splitted for any $\sigma \in \mathbb{N}$.
\end{corollary}

Having in mind Theorem~\ref{Tw1} and Corollary~\ref{Cor1} we say that the
solution $x_{i}(\cdot )$ occurring in the FSS $\bigl\{x_{1}(\cdot ),\dots
,x_{s}(\cdot )\bigr\}$ is \emph{broken away} if it is $\sigma $-broken
away for some $\sigma \in \mathbb{N}$. Accordingly, a FSS is called \emph{%
splitted} if it is $\sigma $-splitted for some $\sigma \in \mathbb{N}$.

\begin{definition}
System~\eqref{2} that has a splitted normal FSS is called  
\emph{a splitted system}.
\end{definition}

The next example shows that there are systems which are not splitted, in
particular such that they have no splitted FSS.

\begin{example}
\textrm{Consider system~\eqref{2} with 
\begin{equation*}
A( n) = 
\begin{pmatrix}
1 & 1 \\ 
0 & 1%
\end{pmatrix}%
,\quad n\in \mathbb{N}.
\end{equation*}%
Obviously,%
\begin{equation*}
X_{A}( n,m) =%
\begin{pmatrix}
1 & n-m \\ 
0 & 1%
\end{pmatrix}%
.
\end{equation*}
Let $x(\cdot)$ be any nontrivial solution of system~\eqref{2} with given $%
A(\cdot)$ and let $x(1)=\mathop{\rm col}(\alpha,\beta)$. Then $%
x(n)=X_A(n,1)x(1)=\mathop{\rm col}(\alpha+\beta n,\beta)$. If $\beta\ne 0$,
then $\sphericalangle\bigl(x(n),e_1\bigr)\doteq\vartheta\to 0$ or $%
\vartheta\to\pi$ as $n\to\infty$, since 
\begin{equation*}
\cos\vartheta=\frac{\langle x(n),e_1\rangle}{\|e_1\|\,\|x(n)\|}=\frac{%
\alpha+\beta n}{\sqrt{\beta^2+(\alpha+\beta n)^2}}\to \frac{\beta}{|\beta|}
\end{equation*}
as $n\to\infty$. On the other hand, if $\beta=0$, then $\vartheta=0$ for all 
$n\in\mathbb{N}$. Hence the angle between any two solutions from any FSS
tends to $0$ or $\pi$. By Definition~\ref{ras4l} it means that any FSS of
this system is not splitted.}
\end{example}

\begin{definition}[{see \protect\cite[p.\thinspace 15]{Halanay},\,%
\protect\cite[p.\thinspace 100]{GAI}}]
\textrm{Let $L(\cdot)$ be a Lyapunov sequence. A linear transformation 
\begin{equation}
y=L(n)x,\quad n\in \mathbb{N},  \label{LyapTr}
\end{equation}%
where $x,y\in\mathbb{R}^{s}$, is called} 
\emph{a Lyapunov transformation.}
\end{definition}

\begin{theorem}
\it
A Lyapunov transformation preserves the property of 
a solution being broken away.
\end{theorem}

\emph{Proof.}
Suppose the solution $x_{i}(\cdot )$ occurring in the FSS $\bigl\{x_{1}(\cdot
),\dots ,x_{s}(\cdot )\bigr\}$ is broken away. Let us apply a Lyapunov
transformation~\eqref{LyapTr} to~\eqref{2}. We shall show that the solutions 
$y_{i}(\cdot )$ from the FSS $\bigl\{y_{1}(\cdot ),\dots ,y_{s}(\cdot )%
\bigr\}
$, where $y_{j}(\cdot )\doteq L(\cdot )x_{j}(\cdot )$, $j=1,\dots ,s$, of
the transformed system, is broken away.

We take an arbitrary $\sigma \in \mathbb{N}$, and by $\psi _{i}(n)$ we
denote the angle between $y_{i}( n) $ and the linear span $%
L(n)V_{i}(n)$ of the vectors $y_{k}(n)$, $k\neq i$, $n\in\mathbb{N}$. It
follows from Lemma~\ref{L1} that%
\begin{equation*}
\psi _{i}(n)\ge \frac{2}{\pi }\varphi _{i}( n) \varkappa
^{1-s}( L(n)) 
=\frac{2}{\pi }\varphi _{i}( n)
\left\| L( n) \right\| ^{1-s}\left\| L^{-1}(
n) \right\| ^{1-s}.
\end{equation*}%
Since $L(\cdot)$ is a Lyapunov sequence, it follows that there exists a $c>0$
such that%
\begin{equation*}
\frac{2}{\pi }\left\| L( n) \right\| ^{1-s}\left\|
L^{-1}( n) \right\| ^{1-s}>c
\end{equation*}%
for all $n\in \mathbb{N}$. Consequently,%
\begin{equation*}
\psi _{i}(n)\ge c\varphi _{i}( n)
\end{equation*}%
for all $n\in \mathbb{N}$. For $\alpha \in \left( 0,\frac{\pi }{2}\right] $,
we set 
\begin{equation*}
L\Gamma _{i}^{\alpha }\doteq\left\{ j\in \mathbb{N}\colon\psi _{i}(j\sigma
)\ge \alpha \right\} ,
\end{equation*}%
\begin{equation*}
L\Gamma _{i}^{\alpha }( k)\doteq L\Gamma _{i}^{\alpha }\cap
\left\{ 1,\dots,k\right\} ,
\end{equation*}%
\begin{equation*}
Lg_{i}^{\alpha }(n)=\frac{1}{k}\sum_{j\in L\Gamma _{i}^{\alpha }(
k) }1.
\end{equation*}%
If $j\in \Gamma _{i}^{\alpha }$, then 
\begin{equation*}
\psi _{i}(j\sigma )\ge c\varphi _{i}(j\sigma )\ge c\alpha ,
\end{equation*}%
i.e., $j\in L\Gamma _{i}^{c\alpha }$. Consequently,%
\begin{equation*}
\Gamma _{i}^{\alpha }( k) \subset L\Gamma _{i}^{c\alpha}(k)
\end{equation*}%
and 
\begin{equation*}
g_{i}^{\alpha }(k)\le Lg_{i}^{c\alpha }(k)
\end{equation*}%
for all $k\in \mathbb{N}$.

Let $\gamma \in \left( 0,\frac{\pi }{2}\right] $ and 
assume the sequence
$( k_{j})_{j\in \mathbb{N}}$ satisfies the relations 
\begin{equation*}
\lim_{j\to \infty }\frac{\ln \|x_{i}(k_{j}\sigma) \|}{k_{j}\sigma }=\lambda[
x_{i}]
\end{equation*}%
and%
\begin{equation*}
\lim_{j\to\infty }g_{i}^{\gamma }(k_{j})>0.
\end{equation*}%
Then 
\begin{equation*}
\lim_{j\to \infty }\frac{\ln \left\| y_{i}( k_{j}\sigma
) \right\| }{k_{j}\sigma }=\lambda [ y_{i}] ,
\qquad
\limsup_{j\to \infty }Lg_{i}^{c\gamma }(k_{j})\ge
\lim_{j\to \infty }g_{i}^{\gamma }(k_{j})>0.
\end{equation*}%
We extract a subsequence $( k_{j_{m}}) _{m\in \mathbb{N}}$ of the
sequence $( k_{j}) _{j\in \mathbb{N}}$ on which the upper limit
$\limsup\limits_{j\to \infty }Lg_{i}^{c\gamma }(k_{j})$ is realized.
Then%
\begin{equation*}
\lim_{m\to \infty }\frac{\ln \left\| y_{i}( k_{j_{m}}\sigma
) \right\| }{k_{j_{m}}\sigma }=\lambda [ y_{i}] 
\end{equation*}%
and%
\begin{equation*}
\lim_{m\to \infty }Lg_{i}^{c\gamma }(k_{j_{m}})>0.
\end{equation*}%
It means that the solution $y_{i}(\cdot )$ is $\sigma $-broken away. The
proof of the theorem is complete.
\hfill $\square$

\begin{corollary}
\it
A Lyapunov transformation preserves the splitting property of a system.
Moreover, a splitted FSS is transformed into a splitted FSS.
\end{corollary}

\begin{remark}
\label{R1} \textrm{Each system~\eqref{2} that can be reduced by a Lyapunov
transformation to a diagonal form is splitted, since the solutions $x_1(\cdot),\dots,x_s(\cdot)$ 
of a diagonal system with the initial conditions $x_i(1)=e_i$, $i=1,\dots,s$, form a normal FSS of this system and 
preserve constant angles between themselves (equal to $\frac{\pi }{2}$).}
\end{remark}

\section{Basic result} \label{sect5}
In this section, we prove the main property of splitted systems.

\begin{theorem}
\it
\label{T6}Suppose that system~\eqref{2} has a splitted FSS 
$x_{1}(\cdot),\dots,x_{s}(\cdot)$. Then there exist $\beta>0$ and $\delta >0$ such
that for any $\xi_{i}\in [-\delta,\delta]$, $i=1,\dots ,s$,  there exists an admissible
multiplicative perturbation $R(\cdot )$, satisfying the estimate%
\begin{equation}
\| R-I\| _{\infty }\le \beta \max 
\bigl\{|\xi _{i}| \colon i=1,\dots ,s\bigr\}  
\label{X7}
\end{equation}%
and such that the solutions $\overline{x}_{i}(\cdot )$, $i=1,\ldots ,s$, of
system~\eqref{M4} with the initial conditions $\overline{x}_{i}(1)=x_{i}(1),$
$i=1,\dots ,s$, satisfy the relations%
\begin{equation*}
\lambda [ \overline{x}_{i}] =\lambda [ x_{i}]+\xi _{i},\quad i=1,\dots ,s.
\end{equation*}%
\end{theorem}

\emph{Proof.}
Fix $\sigma=1$. Since the solutions $x_{i}(\cdot )$, 
$i=1,\dots ,s$, are broken away, it follows that there exist a number 
$\gamma \in \left( 0,\frac{\pi }{2}\right] $ and realizing sequences 
$\bigl(k_{j}(i)\bigr)_{j\in \mathbb{N}}\subset \mathbb{N}$
for solutions $x_{i}(\cdot )$, $i=1,\dots ,s$, such that%
\begin{equation*}
\rho _{i}=\lim_{j\to \infty }g_{i}^{\gamma }( k_{j}(i)) >0
\end{equation*}%
for any $i=1,\dots ,s$. Note that the inequality $\rho _{i}\le 1$ is always
valid, since%
\begin{equation*}
\sup \left\{ g^{\gamma}_{i}(k)\colon k\in \mathbb{N},\ i=1,\dots ,s\right\} \le 1.
\end{equation*}%
This, together with Lemma~\ref{L2}, implies that each function 
\begin{equation*}
\Lambda _{i}^{\gamma }( \mu ) \doteq \limsup_{k\to \infty
}( f_{i}(k)+\mu g_{i}^{\gamma }(k))  
\end{equation*}%
satisfies the estimate%
\begin{equation}
\lambda [ x_{i}] +\mu \ge \Lambda _{i}^{\gamma }( \mu
) \ge \lambda [ x_{i}] +\rho \mu ,
\label{X2}
\end{equation}%
where%
\begin{equation*}
\rho =\min \left\{ \rho _{i}\colon i=1,\dots ,s\right\} ,
\end{equation*}%
and for each $t\ge 0$, there exists a $\mu _{t}^{i}\in [ 0,\rho ^{-1}t%
] $ such that%
\begin{equation}
\Lambda _{i}^{\gamma }( \mu _{t}^{i}) =\lambda [ x_{i}]
+t.  \label{X3}
\end{equation}%
Since $\gamma $ is a 
fixed number throughout the proof, from now on we omit the
superscript~$\gamma $.

Let us fix $r\in ( 0,1) $ and $\delta _{1}\in ( 0,\ln (
L_{1}+1) ) $, where%
\begin{equation*}
L_{1}=\frac{r\sin \gamma }{s}
\end{equation*}%
and denote%
\begin{equation*}
L=\frac{L_{1}}{\delta _{1}}.
\end{equation*}%
It is easy to verify that 
\begin{equation*}
\left| \exp ( \delta _{1}) -1\right| <L_{1},\quad
\left| \exp ( -\delta _{1}) -1\right| <L_{1}
\end{equation*}%
and 
\begin{equation*}
\left| \exp ( \tau ) -1\right| \le L\left| \tau
\right| 
\end{equation*}%
for all $\tau \in ( -\infty ,\delta _{1}] $. 
We set 
\begin{equation*}
\delta =\frac{\delta _{1}\rho }{3}
\end{equation*}%
and take an arbitrary $\xi_{i}\in [-\delta,\delta]$, $i=1,\dots,s$. 
Let $\eta\doteq\min\bigl\{\xi_i\colon i=1,\dots,s\bigr\}$ and $\zeta_i\doteq\xi_i-\eta$. Then
$|\eta|\le\delta$ and
$0=\xi_i-\xi_i\le\xi_i-\eta=\zeta_i\le\xi_i+|\eta|\le 2\delta$,
that is,
\begin{equation*}
0\le\zeta_i\le 2\delta.
\end{equation*}
Let
\begin{equation*}
\varepsilon \doteq\max \bigl\{|\xi_{i}| \colon i=1,\dots ,s\bigr\}.
\end{equation*}%
Then 
\begin{equation*}
|\eta| \le \varepsilon  
\label{X4}
\end{equation*}%
and%
\begin{equation*}
\zeta_i\le|\xi_i|+|\eta|\le 2\varepsilon,\quad i=1,\dots,s.   
\label{X5}
\end{equation*}%
Let the quantities $\mu _{i},$ $i=1,\dots ,s,$ be
obtained from the conditions 
\begin{equation*}
\Lambda _{i}(\mu _{i})=\lambda [ x_{i}] +\zeta _{i}.
\end{equation*}%
Since $\zeta_{i}\ge 0$, it follows that the numbers $\mu _{i}$ are well
defined by~\eqref{X3}, and by~\eqref{X2} we have the estimates 
\begin{equation*}
\mu _{i}\ge \Lambda _{i}(\mu _{i})-\lambda [ x_{i}] =\zeta_{i}\ge \rho \mu _{i}\ge 0
\end{equation*}%
for all $i=1,\dots ,s$.

For each $n\in \mathbb{N}$, we introduce a matrix $R(n)\in \mathbb{R}%
^{s\times s}$ by the formulas%
\begin{equation}
R(n)x_{i}( n) =x_{i}( n) \exp \bigl(s_{i}(
n) \bigr),\quad i=1,\dots ,s,  
\label{N7}
\end{equation}%
where 
\begin{equation*}
s_{i}( n) =
\begin{cases}
\eta +\mu _{i} & \textrm{for\ } n\in \Gamma _{i},\\
\eta & \textrm{for\ } n\not\in \Gamma _{i}.
\end{cases}
\end{equation*}%
Then we have the estimates%
\begin{equation*}
|s_{i}(n)|\le|\eta|+\mu_{i}\le|\eta| +\frac{\zeta_{i}}{\rho}\le
\delta+\frac{2\delta}{\rho}=\frac{\delta}{\rho}(\rho+2)\le\frac{3\delta}{\rho}=
\delta_{1}
\end{equation*}%
for $i=1,\dots ,s$ and $n\in \mathbb{N}$. This, together with the definition
of $\delta_{1}$, implies that 
\begin{equation*}
\bigl|\exp\bigl(s_{i}(n)\bigr)-1\bigr| \le L_{1}
\end{equation*}%
and 
\begin{equation*}
\bigl|\exp\bigl(s_{i}(n)\bigr)-1\bigr|\le
L\bigl| s_{i}( n) \bigr| \le 
L\left(|\eta|+\frac{\zeta_i}{\rho}\right)\le
L\left(|\eta|+\frac{2\varepsilon}{\rho}\right)\le
\left( 1+\frac{2}{\rho }%
\right) L\varepsilon.
\end{equation*}%
By the definition of 
a FSS, the vectors $x_{1}(n),\dots,x_{s}(n)$ are linearly independent for each $n\in \mathbb{N}$
and, by~\eqref{N7}, they are eigenvectors of the matrix $R(n)$. This implies
that $R(n)$ is a matrix of simple structure \cite[p.\,239, Proposition 2]{Lancaster} 
and therefore it can be represented by the sum%
\begin{equation*}
R(n)=\sum\limits_{i=1}^{s}P_{n}^{i}\exp\bigl(s_{i}(n)\bigr),
\end{equation*}%
where the root projections $P_{n}^{i}$ are given by the conditions%
\begin{equation*}
P_{n}^{i}x_{i}( n) =x_{i}( n)
\end{equation*}%
and%
\begin{equation*}
P_{n}^{i}x_{j}( n) =0
\end{equation*}%
for $j\neq i$. Moreover 
\begin{equation*}
\sum\limits_{i=1}^{s}P_{n}^{i}=I
\end{equation*}%
for all $n\in \mathbb{N}$. By Lemma~\ref{L3}, we have the estimate 
\begin{equation*}
\left\| P_{n}^{i}\right\| \le \frac{1}{\sin \gamma }
\end{equation*}%
and hence the inequalities 
\begin{equation}
\begin{split}
\left\| R(n)-I\right\| &=  \Bigl\|
\sum\limits_{i=1}^{s}P_{n}^{i}\bigl(\exp\bigl(s_{i}(n)\bigr)-1\bigr)\Bigr\|\\ 
&\le \sum\limits_{i=1}^{s}\left\|
P_{n}^{i}\right\| \left| \exp \bigl( s_{i}( n) \bigr)-1\right| <s\frac{L_{1}}{\sin \gamma }=r. 
\end{split}
\label{X6}
\end{equation}%
Moreover,%
\begin{equation*}
\left\| R(n)-I\right\| \le s\frac{\left| \exp ( s_i(n)) -1\right| }{\sin \gamma }\le \frac{s}{\sin \gamma }\left( 1+%
\frac{2}{\rho }\right) L\varepsilon .
\end{equation*}%
Hence, for the sequence $R(\cdot )=\bigl( R(n)\bigr)_{n\in\mathbb{N}}$ we
have%
\begin{equation*}
\left\| R-I\right\| _{\infty }=\sup_{n\in \mathbb{N}}\left\|
R(n)-I\right\| \le \beta \varepsilon ,
\end{equation*}%
where%
\begin{equation*}
\beta =\frac{Ls\left( 1+\frac{2}{\rho }\right) }{\sin \gamma }.
\end{equation*}%
This proves~\eqref{X7}.

Moreover, we have%
\begin{equation*}
X_{AR}( n+1,n) =X_{A}( n+1,n) R(n)
\end{equation*}%
for all $n\in \mathbb{N}$. Since $r\in ( 0,1) $, then the
condition $\left\| H-I\right\| <r$ implies that $H$ is invertible and%
\begin{equation*}
\left\| H\right\| \le r+1,
\end{equation*}%
\begin{equation*}
\left\| H^{-1}\right\| \le \frac{1}{1-r}
\end{equation*}%
whatever $H\in \mathbb{R}^{s\times s}$ is given~\cite[p.\,301]{HJ}. 
By~\eqref{X6} this
in turn implies that the sequence $R(\cdot )$ is an admissible
multiplicative perturbation.

Consider the FSS $\overline{x}_{i}(\cdot )$, $i=1,\ldots ,s$, of~\eqref{M4}
with such a perturbation with the initial conditions%
\begin{equation*}
\overline{x}_{i}( 1) =x_{i}( 1) ,\quad i=1,\ldots ,s.
\end{equation*}%
For every natural $i\le s$ and $k\ge 2$ we have the equalities
\begin{equation*}
\begin{split}
\overline{x}_{i}(k)&=X_{AR}(k,1)\overline{x}_{i}(1)=X_{AR}(k,1)x_{i}(1)=\prod_{j=1}^{k-1}X_{AR}(j+1,j)x_{i}(1)\\
&=\prod_{j=1}^{k-1}X_{A}(j+1,j)R(j)x_{i}(1)=
\prod_{j=1}^{k-1}X_{A}(j+1,j)\exp\bigl(s_i(j)\bigr)x_{i}(1)\\
&=\exp\bigl(s_i(1)+\ldots+s_{i}(k-1)\bigr)X_A(k,1)x_i(1)=\exp\Bigl(\sum_{j=1}^{k-1}s_i(j)\Bigr)x_i(k).
\end{split}
\end{equation*}
It follows that
the Lyapunov exponents of these solutions satisfy the relations
\begin{equation*}
\begin{split}
\lambda[\overline{x}_{i}] 
&=\limsup_{k\to \infty }\frac{1}{k}\ln\|\overline{x}_{i}(k)\| 
=\limsup_{k\to \infty }\frac{1}{k}\Bigl(\ln\|x_{i}(k)\| +\sum\limits_{j=1}^{k-1}s_{i}(j)\Bigr)\\
&=\limsup_{k\to\infty}\Bigl(f_{i}(k)+\frac{\mu_{i}N_{i}(k-1)}{k}+\frac{\eta(k-1)}{k}\Bigr)\\ 
&=\limsup_{k\to \infty }\Bigl(f_{i}(k)+\frac{\mu_{i}N_{i}(k)}{k}+\frac{\mu_i(N_i(k-1)-N_i(k))}{k}+\eta-\frac{\eta}{k}\Bigr)\\
&=\limsup_{k\to \infty }\bigl(f_{i}(k)+\mu_{i}g_{i}(k)\bigr) +\eta= 
\Lambda _{i}( \mu _{i}) +\eta =\lambda [ x_{i}] +\zeta_{i}+\eta
=\lambda [ x_{i}] +\xi_{i}
\end{split}
\end{equation*}%
for $i=1,\dots ,s$. 
\hfill $\square$

\section{Applications} \label{sect6}
In this section, we prove several results that follow from Theorem~\ref{T6} and demonstrate the importance of the introduced 
concept of splitted systems for studying the behavior of the Lyapunov spectrum  under the action of 
small perturbations. 

\begin{theorem}
\label{T7}
\it
If system~\eqref{2} has a splitted FSS which is not normal, then
the Lyapunov spectrum of system~\eqref{2} is not stable.
\end{theorem}

\emph{Proof.}
Let $\bigl\{x_{1}(\cdot ),\dots ,x_{s}(\cdot )\bigr\}$ be a splitted FSS of system~\eqref{2}
which is not normal. 
Denote 
\begin{equation*}
\mu _{i}=\lambda [ x_{i}] ,\quad i=1,\dots ,s.
\end{equation*}%
Without loss of generality, we assume that the sequence $\mu =\bigl(\mu
_{1},\dots ,\mu _{s}\bigr)\in \mathbb{R}_{\leqslant }^{s}$. Note that among
the numbers $\mu _{1},\dots ,\mu _{s}$, there are equal ones, because
otherwise the FSS $\bigl\{x_{1}(\cdot ),\dots ,x_{s}(\cdot )\bigr\}$ would
be incompressible, and therefore normal. If 
\begin{equation*}
\lambda (A)=\bigl(\lambda _{1}(A),\ldots ,\lambda _{s}(A)\bigr)\in \mathbb{R}%
_{\leqslant }^{s}
\end{equation*}%
is the Lyapunov spectrum of~\eqref{2}, then $\mu \neq \lambda (A)$.

Choose a number $\alpha >0$ so small that the sets 
${\mathcal{O}}_{\alpha}(\mu )$ and ${\mathcal{O}}_{\alpha }\bigl(\lambda (A)\bigr)$ do not
intersect. Take an arbitrary positive $\varepsilon <\min \{\alpha ,\delta \}$, 
where $\delta $ is from Theorem~\ref{T6}. With the selected $\varepsilon$,
in the neighborhood of ${\mathcal{O}}_{\varepsilon }(\mu )$ there is a
sequence of numbers $\mu ^{\prime }=\bigl(\mu _{1}^{\prime },\dots ,\mu_{s}^{\prime }\bigr)$, 
all of whose elements are different, i.e., $\mu^{\prime }\in \mathbb{R}_{<}^{s}$. 
Take 
$\xi_i\doteq\mu'_i-\mu_i$, $i=1,\dots,s$. 
Then $|\xi_i|<\varepsilon<\delta$, therefore, by Theorem~\ref{T6}, there exist an admissible multiplicative perturbation 
$R(\cdot)$ satisfying the estimate 
$\|R-I\|_{\infty}<\beta\varepsilon$ and such that the system~\eqref{M4} with
such a perturbation has a FSS $\bigl\{\overline{x}_{1}(\cdot),\dots,\overline{x}_{s}(\cdot)\bigr\}$ 
such that 
\begin{equation*}
\lambda[\overline{x}_{i}]=\lambda[{x}_{i}]+\xi_i=\mu_i+\xi_i=\mu'_{i},\quad i=1,\dots,s.
\end{equation*}
 Since the numbers $\mu _{i}^{\prime }$
are pairwise distinct, this FSS is incompressible and therefore normal.
Hence, the set $\mu ^{\prime }$ is the Lyapunov spectrum of system~\eqref{M4}. 
Obviously, $\mu ^{\prime }\not\in {\mathcal{O}}_{\alpha }\bigl(\lambda (A)\bigr)$ 
for all such $\varepsilon $, which means that the Lyapunov spectrum
of system~\eqref{2} is not stable.
\hfill $\square$

\begin{corollary}
\label{sled}
\it
If the Lyapunov spectrum of system~\eqref{2} is stable, then each splitted
FSS is normal.
\end{corollary}

It can be easily seen that stability of Lyapunov spectrum
is equivalent to continuity of the map $R(\cdot)\mapsto\lambda(AR)$ at the point $R(n)\equiv I$, $n\in\mathbb N$.
Theorem~\ref{T6} can be also used  to study some other  properties of this map.

\begin{definition}
\label{open}
The Lyapunov spectrum of system~\eqref{2} is called \emph{open},
if the mapping $\lambda(AR)\colon\mathcal R\to\mathbb R^{s}_{\le}$ is open at the point $R(n)\equiv I$, $n\in\mathbb N$,
that is, for any $\varepsilon>0$, there exists $\gamma=\gamma(\varepsilon)>0$ such that the inclusion
\begin{equation*}
\mathcal O_{\gamma}\bigl(\lambda(A)\bigr)\subset\lambda\bigl(\mathcal R_{\varepsilon}\bigr)
\end{equation*}
holds.
\end{definition}

\begin{remark}
It was proved in~\cite[Theorem~3]{BP} that if the Lyapunov spectrum
of system~\eqref{2} is stable, then it is open.
Here we obtain another sufficient condition for 
the Lyapunov spectrum to be open, expressed in terms
of the splitness of system~\eqref{2}.
\end{remark}

\begin{theorem}
\label{newtheor2}
\it
If system~\eqref{2} is splitted and has a non-multiple Lyapunov spectrum, i.e.
$\lambda(A)\in\mathbb R^{s}_{<}$, then the Lyapunov spectrum of this system is open.
\end{theorem}

\emph{Proof.}
Let $\lambda(A)=\bigl(\lambda_1,\dots,\lambda_s\bigr)\in\mathbb R^{s}_{<}$ and
$\bigl\{x_1(\cdot),\dots,x_s(\cdot)\bigr\}$ be a normal splitted FSS of system~\eqref{2}, such that $\lambda[x_i]=\lambda_j$, $i=1,\dots,s$.
Denote
\begin{equation*}
\eta\doteq\frac{1}{3}\min\bigl\{\lambda_{i+1}-\lambda_i\colon i=1,\dots,s-1\bigr\}.
\end{equation*}
For arbitrary $\varepsilon>0$, we put
\begin{equation*}
\gamma=\gamma(\varepsilon)=\min\{\eta,\, \varepsilon/\beta,\, \delta/\beta,\, \delta\},
\end{equation*}
where $\delta>0$ and
$\beta>0$ are the quantities from Theorem~\ref{T6}.
Take any $\mu=\bigl(\mu_1,\dots,\mu_s\bigr)\in\mathcal O_{\gamma}\bigl(\lambda(A)\bigr)$ and prove that
$\mu\in\lambda\bigl(\mathcal R_{\varepsilon}\bigr)$.
Let $\xi_i\doteq\mu_i-\lambda_i$.
Then $|\xi_i|<\gamma\le\delta$.
By theorem~\ref{T6} there exists
an admissible multiplicative perturbation $R(\cdot)$
that ensures the equality
\begin{equation*}
\lambda[\overline{x}_{i}] =\lambda[x_{i}]+\xi_{i}=\lambda_i+\xi_{i}=\mu_i
\end{equation*}%
for the solution of system~\eqref{M4} with the initial condition $\overline{x}_{i}(1)=x_{i}(1)$,
and such that the inequality
\begin{equation*}
\|R-I\|_{\infty}<\beta\max\bigl\{|\xi_i|\colon i=1,\dots,s\bigr\}<\beta\gamma\le\varepsilon
\end{equation*}
holds, i.e.
$R(\cdot)\in\mathcal R_{\varepsilon}$.

Consider the FSS $\bigl\{\overline x_1(\cdot),\dots,\overline x_s(\cdot)\bigr\}$ of system~\eqref{M4}.
Let us note that
\begin{equation*}
\begin{split}
&\lambda[\overline x_{i+1}]-\lambda[\overline x_{i}]
=\lambda_{i+1}-\lambda_{i}+\xi_{i+1}-\xi_{i}\\
&\ge 3\eta-|\xi_{i+1}|-|\xi_{i}|\ge 3\eta-2\delta\ge 3\eta-2\eta=\eta>0
\end{split}
\end{equation*}
for $i\in\{1,\dots,s-1\}$.
Hence the numbers $\lambda[\overline x_{1}],\dots,\lambda[\overline x_{s}]$ are pairwise different,
so the FSS $\bigl\{\overline x_1(\cdot),\dots,\overline x_s(\cdot)\bigr\}$ of system~\eqref{M4} is normal and
\begin{equation*}
\lambda(AR)=\bigl(\lambda[\overline x_{1}],\dots,\lambda[\overline x_{s}]\bigr)=\mu.
\end{equation*}
It means that $\mu\in\lambda\bigl(\mathcal R_{\varepsilon}\bigr)$.
\hfill $\square$

The property of openness of the Lyapunov spectrum of system~\eqref{2} can be interpreted as
the property of local assignability of the Lyapunov spectrum of system~\eqref{M4} under
the action of a multiplicative perturbation $R(\cdot)$,
which in this context is considered as a matrix control.

\begin{definition}
\label{lokupr}
The Lyapunov spectrum of system~\eqref{M4} is called \emph{locally assignable} if for any $\varepsilon>0$ there exists such
a $\gamma=\gamma(\varepsilon)>0$ that for any $\mu\in\mathcal O_{\gamma}\bigl(\lambda(A)\bigr)$ there is
a matrix control $R(\cdot)\in\mathcal{R}_{\varepsilon}$ such that $\lambda(AR)=\mu$.
\end{definition}

It is clear that the property of openness of the Lyapunov spectrum of system~\eqref{2}
coincides with the property of local assignability of the Lyapunov spectrum of system~\eqref{M4}; therefore, the following corollary holds.

\begin{corollary}
\it
If system~\eqref{2} is splitted and has a non-multiple Lyapunov spectrum,
then the Lyapunov spectrum of system~\eqref{M4} is locally assignable.
\end{corollary}

\section{Examples} \label{sect7}

In this section we shall present examples that show that there are both
systems with broken away normal FSS and systems having broken away FSS which
is not normal. We shall also show that the assumptions of the proved
theorems can be effectively verified.

In the first example we present a system with splitted FSS which is not
normal and therefore its Lyapunov spectrum is not stable. In this example we
use the following result.

\begin{lemma}
\it
We have%
\begin{equation}
\limsup_{n\to \infty }\sin \ln n=1,
\label{sin1}
\end{equation}%
\begin{equation}
\liminf_{n\to \infty }\sin \ln n=-1. 
\label{sin2}
\end{equation}
\end{lemma}

\emph{Proof.}
Consider the following sequences 
\begin{equation*}
t_{k}=\exp \left( 2k+\frac{1}{2}\right) \pi ,\ n_{k}=[ t_{k}]
,\quad k\in \mathbb{N}.
\end{equation*}%
We have%
\begin{equation*}
\begin{split}
1\ge & \limsup_{n\to \infty }\sin \ln n\ge \limsup_{k\to
\infty }\sin \ln n_{k}
=\limsup_{k\to \infty }\sin \bigl(\ln t_{k}+\ln (n_{k}/t_{k})\bigr) \\
= & \limsup_{k\to \infty }\sin \Bigl(\sin \bigl(\ln t_{k}\bigr)\cos %
\bigl(\ln (n_{k}/t_{k})\bigr) 
+\cos \bigl(\ln t_{k}\bigr)\sin \bigl(\ln
(n_{k}/t_{k})\bigr)\Bigr) \\
= &\limsup_{k\to \infty }\cos \bigl(\ln (n_{k}/t_{k})\bigr)=\cos \ln %
\bigl(\lim_{k\to \infty }(n_{k}/t_{k})\bigr) 
=\cos \ln 1=1.
\end{split}
\end{equation*}%
Considering sequences 
\begin{equation*}
t_{k}=\exp \left( 2k+\frac{3}{2}\right) \pi ,\text{ }n_{k}=[ t_{k}%
] ,\text{ }k\in \mathbb{N}\text{,}
\end{equation*}%
we may prove in a similar way the equality~\eqref{sin2}.
\hfill $\square$

\begin{example}
\label{ex2} 
\textrm{Let us consider system~\eqref{2} with $s=2$ and a
diagonal matrix 
\begin{equation*}
A(n)=\mathop{\rm diag}\bigl(a_1(n),a_2(n)\bigr),\quad n\in\mathbb{N},
\end{equation*}
where 
\begin{align*}
a_1(n) & =  \exp\bigl(n\sin\ln n-(n+1)\sin\ln(n+1)\bigr), \\
a_2(n) & =  \exp\Bigl(2\bigl((n+1)\sin\ln(n+1)-n\sin\ln n\bigr)\Bigr).
\end{align*}
}

\textrm{The sequence $A(\cdot )$ is a Lyapunov sequence since%
\begin{equation*}
\bigl|(n+1)\sin \ln (n+1)-n\sin \ln n\bigr|\le \sqrt{2},\quad n\in \mathbb{N%
}.
\end{equation*}%
To obtain the last inequality it is enough to apply the Lagrange mean value
theorem to the function $f\colon [ n,n+1] \to \mathbb{R}$%
, $f(t)=t\sin \ln t$. It is easy to see that the matrix%
\begin{equation*}
\Phi(n,1)=\left( 
\begin{array}{cc}
\exp \bigl(-n\sin \ln n\bigr) & 0 \\ 
0 & \exp \bigl(2n\sin \ln n\bigr)%
\end{array}%
\right) ,\quad n\in \mathbb{N},
\end{equation*}%
is a fundamental matrix of this system. By Remark~\ref{R1} the corresponding
FSS is normal. By~\eqref{sin1},\thinspace \eqref{sin2} we have%
\begin{equation*}
\lambda (A)=(1,2).
\end{equation*}%
Consider now the FSS $\bigl\{x_{1}(\cdot ),x_{2}(\cdot )\bigr\}$, where 
\begin{equation*}
x_{1}(n)=\left( 
\begin{array}{c}
\exp \bigl(-n\sin \ln n\bigr) \\ 
\exp \bigl(2n\sin \ln n\bigr)%
\end{array}%
\right) ,
\quad
x_{2}( n) =\left( 
\begin{array}{c}
0 \\ 
\exp \bigl(2n\sin \ln n\bigr)%
\end{array}%
\right) ,\quad n\in \mathbb{N}.
\end{equation*}%
From~\eqref{sin1} and~\eqref{sin2} it follows that $\lambda [ x_{1}%
] =\lambda [ x_{2}] =2$ and that 
\begin{equation*}
n_{k}=[ t_{k}] +1,\quad k\in \mathbb{N},
\end{equation*}%
is a realizing sequence for $x_{1}(\cdot )$ and $x_{2}(\cdot )$, where $%
t_{k}=\exp \left( 2k+\frac{1}{2}\right) \pi .$ Therefore FSS $\bigl\{%
x_{1}(\cdot ),x_{2}(\cdot )\bigr\}$ is not normal. We shall show that it is
splitted. }

\textrm{Denote by $\varphi _{1}(n)$ the angle between $x_{1}( n) $
and $x_{2}( n) $. After some simple calculations we have 
\begin{equation*}
\cos \varphi _{1}(n)=\bigl(1+\exp (-6n\sin \ln n)\bigr)^{-1/2}.
\end{equation*}%
Let us fix a 
\begin{equation}
c\in ( \sqrt{2}/2,1) .  \label{22.2}
\end{equation}%
Let $\gamma \doteq \arccos c$, then $\gamma \in (0,\pi /4)$. Notice that $%
n\in \Gamma _{1}^{\gamma }(M;1)$ for some $M\in \mathbb{N}$ if and only if $%
n\in \mathbb{N}$, $n\le M$ and 
\begin{equation*}
\cos ^{2}\varphi _{1}(n)\le c^{2},
\end{equation*}%
i.e.,%
\begin{equation*}
1+\exp \bigl(-6n\sin \ln n\bigr)\ge 1/c^{2}.
\end{equation*}%
The last inequality is equivalent to the inequality%
\begin{equation}
6n\sin \ln n\le \ln \left( \frac{c^{2}}{1-c^{2}}\right) .  \label{22.3}
\end{equation}%
By the choice of $c$ we know that $\frac{c^{2}}{1-c^{2}}>1$ and therefore 
\begin{equation*}
\ln \left( \frac{c^{2}}{1-c^{2}}\right) >0.
\end{equation*}%
The last inequality means that each $n\in \mathbb{N}$ satisfying $\sin \ln
n\le 0$ 
also satisfies inequality~\eqref{22.3}
 and therefore 
\begin{equation*}
\left\{ n\in \mathbb{N}\colon n\le M,\sin \ln n\le 0\right\} \subset
\Gamma _{1}^{\gamma }( M;1)
\end{equation*}%
and%
\begin{equation*}
\mathop{\rm card}\left\{ n\in \mathbb{N}\colon n\le M,\sin \ln n\le
0\right\} \le N_{1}^{\gamma }( M;1) ,
\end{equation*}%
where $\mathop{\rm card}B$ denotes the number of elements of the set $B$.
Let us fix $k\in \mathbb{N}$, $k>1$. Then for each 
\begin{equation*}
n\in \bigl[\exp \bigl((2k-1)\pi \bigr),\exp \bigl(2k\pi \bigr)\bigr]\cap 
\mathbb{N}
\end{equation*}%
we have $\sin \ln n\le 0$ and $1<n<n_{k}$, therefore 
\begin{equation*}
\bigl[\exp \bigl((2k-1)\pi \bigr),\exp \bigl(2k\pi \bigr)\bigr]\cap \mathbb{N%
}\subset \Gamma _{1}^{\gamma }( n_{k};1)
\end{equation*}%
and}%
\begin{equation*}
N_{1}^{\gamma }( n_{k};1) \ge \exp ( 2k\pi ) -\exp
( ( 2k-1) \pi ) -1.
\end{equation*}%
Since $n_{k}\le \exp \left( \left( 2k+\frac{1}{2}\right) \pi \right) +1,$
we have 
\begin{equation*}
\limsup_{k\to \infty }\frac{N_{1}^{\gamma }( n_{k};1) }{%
n_{k}}\ge \frac{1-\exp ( -\pi ) }{\exp ( \pi /2) }>0,
\end{equation*}%
which completes the proof. 
\end{example}

In the next example we present a system with 
a stable Lyapunov spectrum such
that each 
of its FSS, which is not normal, is not splitted.

\begin{example}
\label{primer3} Let us consider the system 
\begin{equation}
x(n+1)=\mathop{\rm diag}\bigl(1,2\bigr)x(n),\quad n\in \mathbb{N},\ x\in 
\mathbb{R}^{2}.  \label{pr3}
\end{equation}%
Since the system is time-invariant, it follows from~\cite{BP} that its
Lyapunov spectrum is stable. Note that the normal FM of system~\eqref{pr3}
has the form 
\begin{equation*}
\Phi (n)=\left( 
\begin{array}{cc}
1 & 0 \\ 
0 & 2^{n}%
\end{array}%
\right) ,
\end{equation*}%
and therefore the corresponding normal FSS of system~\eqref{pr3} is
splitted. Let us prove that every FSS of system~\eqref{pr3}, which is not
normal, is not splitted. Let us divide the set of all nontrivial solutions
of system~\eqref{pr3} into two groups. We include 
in the first group all
the solutions $x(\cdot )$ of this system with the initial conditions $%
x(1)=\alpha e_{1}$, where $\alpha \neq 0 $. Such solutions are constant and
their Lyapunov exponents are equal to $0$. 
In the second group we include
all the solutions $x(\cdot )$ of system \eqref{pr3} with initial conditions $%
x(1)=\alpha e_{1}+\beta e_{2}$, where $\beta \neq 0$. They have the form $%
x(n)=\mathop{\rm
col}\bigl(\alpha ,2^{n-1}\beta \bigr)$ and their Lyapunov exponents are
equal to $\ln 2$. Each normal FSS contains solutions from both of these
groups. If FSS is not normal, then it should contain two solutions only from
the second group. We cannot construct a normal FSS from the solutions from
the first group, since any two solutions from the first group are obtained
one from the other by multiplying by some constant, i.e., they are linearly
dependent.

Take any solution $x(\cdot )$ from the second group and calculate the angle
between 
$x(\cdot )$  and the vector $e_{2}$. We have 
\begin{equation*}
\cos \sphericalangle \bigl(e_{2},x(n)\bigr)=\frac{\langle e_{2},x(n)\rangle 
}{\| e_{2}\| \,\| x(n)\| }=
\frac{2^{n-1}\beta }{\sqrt{\alpha ^{2}+4^{n-1}\beta ^{2}}}\to 
\frac{\beta }{|\beta |}\quad \text{for\ }n\to \infty ,
\end{equation*}%
i.e., $\sphericalangle \bigl(e_{2},x(n)\bigr)$ tends to $0$ or $\pi $ when $%
n\to \infty $. It follows that the angle between any two vectors $%
x_{1}(n)$ and $x_{2}(n)$ corresponding to some solutions from the second
group tends to $0$ or $\pi $ when $n\to \infty $. By Definition~\ref{ras4l} it
means that FSS $\bigl\{x_{1}(\cdot ),x_{2}(\cdot )\bigr\}$ which is not
normal is also not splitted.
\end{example}

\section{Conclusions}

In this paper we study the 
stability and openness problems  of Lyapunov spectra
of discrete time-varying linear systems. To investigate this problems we
introduce the concept of broken away solutions and the concept of splitted
systems. We demonstrate some properties of these concepts and 
then applied these properties to the stability and openness problems
  of Lyapunov spectra. One of the main results states
that if the Lyapunov spectrum is stable then 
each splitted fundamental system of solutions is normal. It is 
worth mentioning that this condition does not require a reduction of 
the given system to any special form, but is rather expressed in terms of the system itself. 
Another important result is that if a given system is splitted and has a non-multiple Lyapunov spectrum, 
then it has an open Lyapunov spectrum. 
We expect that the proposed concepts of broken away solutions and splitted
systems may be useful 
in the investigation of other problems 
in the theory of
discrete time-varying linear systems, such as the problem of assignability of
the Lyapunov spectrum. This will be a subject of our further research.

\section*{Acknowledgments} 
The research of the first and third authors was funded by 
the National Science Centre in Poland granted according to decisions DEC-2017/25/B/ ST7/02888 and DEC-2015/19/D/ST7/03679, respectively. 
The research of the fourth author was partially funded by RFBR (project number 20--01--00293) 
and 
by the Ministry of Science and Higher Education of the Russian
Federation in the framework of the state assignment No. 075-00928-21-01 
(project FEWS-2020-0010). 
The work of the fifth author was funded by the Polish National Agency for Academic Exchange  NAWA (the Ulam program) 
granted according to the decision No.~PPN/ULM/2019/1/00287/DEC/1.



\end{document}